# A counterexample to the geometric Chevalley-Warning conjecture

June Huh

### Abstract

We construct a quartic threefold with $\mathbb{L}$-rational singularities which has torsion in its middle homology group. This answers a question of Brown and Schnetz for all fields of characteristic zero.

## 1. Introduction

Let $\mathbb{F}_q$ be the finite field with $q$ elements of characteristic $p$, and let $h$ be a degree $d$ homogeneous polynomial in $\mathbb{F}_q[z_0, \ldots, z_n]$. If $d \leqslant n$, then the classical Chevalley-Warning theorem says that

$$\#(V(h))_{\mathbb{F}_q} = 1 \mod p,$$

where $\#(V(h))_{\mathbb{F}_q}$ is the number of $\mathbb{F}_q$-rational points in the projective hypersurface $V(h)$ defined by $h$. Ax later proved the following refinement [Ax64].

THEOREM 1. *If $h$ is a homogeneous polynomial in $\mathbb{F}_q[z_0, \ldots, z_n]$ of degree $d \leqslant n$, then*

$$\#(V(h))_{\mathbb{F}_q} = 1 \mod q.$$

A conjectural generalization of the Chevalley-Warning theorem was proposed in [BS12] at the level of Grothendieck ring of varieties. Let $k$ be a field. The *Grothendieck ring of $k$-varieties* $K_0(\text{Var}_k)$ is the abelian group generated by classes of finite type schemes over $k$ up to isomorphism, with the relation

$$[X] = [X \backslash Y] + [Y]$$

if $Y \subseteq X$ is a closed subscheme of $X$. The ring structure is defined by the product

$$[X] \cdot [Y] = [X \times_k Y].$$

We denote by 1 the multiplicative identity $[\text{Spec}(k)]$ and by $\mathbb{L}$ the class of the affine line $\mathbb{A}_k^1$.

Brown and Schnetz asked which fields $k$ have the following property [BS12, Question 26].

CHEVALLEY-WARNING PROPERTY.
*If $h$ is a homogeneous polynomial in $k[z_0, \ldots, z_n]$ of degree $d \leqslant n$, then*

$$[V(h)] = 1 \mod \mathbb{L}$$

*in the Grothendieck ring of $k$-varieties.*

In [Lia13] the statement that algebraically closed fields of characteristic zero have the Chevalley-Warning property was given the name "geometric Chevalley-Warning conjecture".

The author was partially supported by NSF grant DMS-0943832.

An interest in the characteristic zero case of Chevalley-Warning comes from its connection with birational geometry. Recall that two $k$-varieties $X$ and $Y$ are said to be *stably birational* if there exist nonnegative integers $m$ and $n$ such that $X \times_k \mathbb{P}^m$ is birational to $Y \times_k \mathbb{P}^n$. When $k$ is an algebraically closed field of characteristic zero, Larsen and Lunts obtained an isomorphism of rings

$$K_0(\text{Var}_k)/(\mathbb{L}) \longrightarrow \mathbb{Z}[\text{SB}_k],$$

where $\text{SB}_k$ is the monoid of stably birational equivalence classes of smooth complete varieties over $k$ [LL03]. This isomorphism maps the class of a smooth complete variety to its stable birational equivalence class. In particular, a smooth complete variety $X$ is *stably rational* (stably birational to a point) if and only if

$$[X] = 1 \mod \mathbb{L}$$

in the Grothendieck ring of $k$-varieties. An example of a smooth projective variety which is not rational but stably rational over $\mathbb{C}$ can be found in [BCSS85].

Over an arbitrary field $k$, Esnault asked the following related question [Bil12, Question 3.7].

QUESTION 2. *Let $h$ be a homogeneous polynomial in $k[z_0, \ldots, z_n]$ of degree $d \leqslant n$. Is it true that $V(h)$ has a $k$-rational point if and only if*

$$[V(h)] = 1 \mod \mathbb{L}$$

*in the Grothendieck ring of $k$-varieties?*

A number of positive results were obtained in [Bil13, Lia13]. On the other hand, Nguyen gave a negative answer to Question 2 for some quasi-algebraically closed and non-algebraically closed fields of characteristic zero [Ngu12], answering the question of Brown and Schnetz for those fields.

We construct a homogeneous polynomial with integer coefficients which shows that no field of characteristic zero has the Chevalley-Warning property.

THEOREM 3. *There is a quartic homogeneous polynomial $h$ in $\mathbb{Z}[z_0, \ldots, z_4]$ such that*

$$[V(h)] \neq 1 \mod \mathbb{L}$$

*in the Grothendieck ring of $k$-varieties for any field $k$ of characteristic zero.*

An explicit construction of $h$ can be found in Section 2.1. The quartic threefold $V(h)$, when viewed over $\mathbb{C}$, has a number of other notable properties:

(i) $V(h)$ has infinitely many $\mathbb{Q}$-rational points (Proposition 5).
(ii) $V(h)$ has $\mathbb{L}$-rational singularities (Definition 6).
(iii) $V(h)$ has torsion in its middle homology group (Proposition 10).
(iv) $V(h)$ is not stably rational (Corollary 11).
(v) $V(h)$ is rationally connected (Remark 13).
(vi) $V(h)$ is a Fano threefold, not stably birational to any smooth Fano threefold (Remark 12).

Furthermore,

(vii) $V(h)$ is simply connected [Laz04, Theorem 3.1.21],
(viii) $V(h)$ has Picard number 1 [Laz04, Example 3.1.25], and
(ix) any smooth model of $V(h)$ has negative Kodaira dimension [CM87, Lemma 1.11].

The author does not know whether this $h$ defines a projective hypersurface with

$$[V(h)] \neq 1 \mod \mathbb{L}$$

in the Grothendieck ring of varieties over a field of positive characteristic.

An approach to proving that smooth quartic threefolds are not stably rational over $\mathbb{C}$ has been recently proposed by Karzhemanov [Kar], but the problem remains open.

### Acknowledgements

The author thanks Kristian Ranestad and Bernd Sturmfels for useful comments on quartic symmetroids. He thanks Igor Dolgachev, Mircea Mustaţă, and Sam Payne for helpful discussions. He thanks Francis Brown and Burt Totaro for valuable comments on previous versions of this paper.

## 2. The quartic threefold $V$

We write $V$ for the quartic threefold $V(h)$. We work over an algebraically closed field $k$ of characteristic zero, unless otherwise stated.

### 2.1 Construction

The quartic threefold $V$ will contain a *quartic symmetroid* as a hyperplane section. By a quartic symmetroid, we mean a quartic surface in $\mathbb{P}^3$ defined by the symmetric determinant

$$\det \begin{pmatrix} l_{11} & l_{12} & l_{13} & l_{14} \\ l_{12} & l_{22} & l_{23} & l_{24} \\ l_{13} & l_{23} & l_{33} & l_{34} \\ l_{14} & l_{24} & l_{34} & l_{44} \end{pmatrix} = 0,$$

where $l_{ij}$ are linear forms in $k[z_0, z_1, z_2, z_3]$. These surfaces, introduced by Cayley, occupy a central place in his memoirs on quartic surfaces [Cay69a, Cay69b, Cay69c]. We summarize the needed properties of quartic symmetroids:

(i) If the linear forms $l_{ij}$ are sufficiently general, then the corresponding quartic symmetroid has 10 isolated nodes.

(ii) Let $\varphi$ be the linear projection from any one of the nodes. Then $\varphi$ is a double covering of $\mathbb{P}^2$ by the quartic symmetroid, ramified along a sextic plane curve.

(iii) The sextic plane curve is the union of two cubic curves $E_1 \cup E_2 \subseteq \mathbb{P}^2$.

For a discussion of the determinantal identity behind the decomposition $E_1 \cup E_2 \subseteq \mathbb{P}^2$, see [Jes16, Chapter IX]. If the linear forms $l_{ij}$ are sufficiently general, then

(iv) the cubic curves $E_1$ and $E_2$ are smooth and meet transversely at 9 points,

(v) there is a smooth conic curve $A$ tangent to each $E_i$ at 3 points, and

(vi) the conic curve $A$ does not intersect $E_1 \cap E_2$.

Quartic symmetroids form a 24-dimensional family in the 34-dimensional space of quartic surfaces. A modern treatment of the above properties of quartic symmetroids can be found in [Cos83].

Choose a quartic symmetroid defined over $\mathbb{Q}$, satisfying the six conditions mentioned above, so that the center of the linear projection $\varphi$ satisfies
$$z_0 = z_1 = z_2 = 0.$$
Then the defining equation of the quartic symmetroid is of the form
$$\alpha(z_0, z_1, z_2)z_3^2 + \beta(z_0, z_1, z_2)z_3 + \gamma(z_0, z_1, z_2) = 0,$$
where $\alpha$, $\beta$, $\gamma$ are defined over $\mathbb{Q}$. The discriminant of the above equation factors over $\mathbb{Q}(i)$
$$\beta(z_0, z_1, z_2)^2 - 4\alpha(z_0, z_1, z_2) \cdot \gamma(z_0, z_1, z_2) = \epsilon_1(z_0, z_1, z_2) \cdot \epsilon_2(z_0, z_1, z_2),$$
and the factors of the discriminant define elliptic curves
$$E_1 = \{\epsilon_1(z_0, z_1, z_2) = 0\} \subseteq \mathbb{P}^2 \quad \text{and} \quad E_2 = \{\epsilon_2(z_0, z_1, z_2) = 0\} \subseteq \mathbb{P}^2.$$
The leading coefficient $\alpha$ defines the smooth conic
$$A = \{\alpha(z_0, z_1, z_2) = 0\} \subseteq \mathbb{P}^2.$$
We choose another sufficiently general quadratic form $\delta(z_0, z_1, z_2)$ over $\mathbb{Q}$, and define
$$D = \{\delta(z_0, z_1, z_2) = 0\} \subseteq \mathbb{P}^2.$$
The conic curve $D$ is smooth and intersects $E_1$, $E_2$, and $A$ transversely. Furthermore, $D \cap E_1 \cap E_2 = \emptyset$.

DEFINITION 4. Let $V \subseteq \mathbb{P}^4$ be the quartic threefold defined by the homogeneous polynomial
$$h := \alpha(z_0, z_1, z_2) \cdot z_3^2 + \beta(z_0, z_1, z_2) \cdot z_3 + \gamma(z_0, z_1, z_2) + \delta(z_0, z_1, z_2) \cdot z_4^2.$$

The notations used in the expression for $h$ are chosen in order to emphasize the similarity with Artin and Mumford's example of a unirational variety which is not rational [AM72, Section 2].

PROPOSITION 5.
 (i) $V$ is singular along the line $L := \{z_0 = z_1 = z_2 = 0\} \subseteq \mathbb{P}^4$,
 (ii) $V$ has 9 isolated nodes ($A_1$-singularities) not contained in $L$, and
(iii) $V$ has no other singularities.

*Proof.* Choose a singular point $p$ of $V$ outside $L$, say with nonzero $z_0$. Introduce affine coordinates
$$x_1 := \frac{z_1}{z_0}, \quad x_2 := \frac{z_2}{z_0}, \quad x_3 := \frac{z_3}{z_0}, \quad x_4 := \frac{z_4}{z_0}$$
for the affine chart $\{z_0 \neq 0\} \subseteq \mathbb{P}^4$. Note that the singular point satisfies
$$\frac{\partial h}{\partial x_3} = \alpha \cdot 2x_3 + \beta = 0 \quad \text{and} \quad \frac{\partial h}{\partial x_4} = \delta \cdot 2x_4 = 0.$$
We claim that $\alpha(p) \neq 0$. If otherwise, $\beta(p) = \gamma(p) = 0$, and hence
$$\beta(p)^2 - 4\alpha(p)\gamma(p) = \epsilon_1(p)\epsilon_2(p) = 0.$$
Since $E_1$ and $E_2$ are smooth, the preceding equation implies that
$$\epsilon_1(p) = \epsilon_2(p) = 0.$$
This contradicts that $A$ does not intersect $E_1 \cap E_2$. Therefore, $\alpha(p) \neq 0$.

Now write the local equation of $V$ at $p$ by

$$\alpha\left(x_3 + \frac{\beta}{2\alpha}\right)^2 - \left(\frac{\beta^2 - 4\alpha\gamma}{4\alpha}\right) + \delta x_4^2 = \alpha\left(x_3 + \frac{\beta}{2\alpha}\right)^2 - \left(\frac{\epsilon_1\epsilon_2}{4\alpha}\right) + \delta x_4^2 = 0.$$

Note that $\epsilon_1(p)\epsilon_2(p) = 0$. We show that $\epsilon_1(p) = \epsilon_2(p) = 0$.

1. If $\epsilon_1(p) = 0$, $\epsilon_2(p) \neq 0$, $\delta(p) = 0$, then take as local parameters at the singular point

$$x_3 + \frac{\beta}{2\alpha}, \quad \epsilon_1, \quad \text{and} \quad \delta.$$

Then the local equation shows that the multiplicity of $V$ at $p$ is 1, a contradiction.

2. If $\epsilon_1(p) = 0$, $\epsilon_2(p) \neq 0$, $\delta(p) \neq 0$, then take as local parameters at the singular point

$$x_3 + \frac{\beta}{2\alpha}, \quad \epsilon_1, \quad \text{and} \quad x_4.$$

Then the local equation shows that the multiplicity of $V$ at $p$ is 1, a contradiction.

Here we used that $E_1$ and $D$ are smooth curves intersecting transversely.

Recall that $D$ does not intersect $E_1 \cap E_2$. Since $\epsilon_1(p) = \epsilon_2(p) = 0$, we have

$$\delta(p) \neq 0 \quad \text{and} \quad z_4(p) = 0.$$

Therefore, the singular points of $V$ outside $L$ are precisely the 9 points satisfying

$$\alpha(z_0, z_1, z_2) \cdot 2z_3 + \beta(z_0, z_1, z_2) = 0, \quad \epsilon_1(z_0, z_1, z_2) = 0, \quad \epsilon_2(z_0, z_1, z_2) = 0, \quad z_4 = 0.$$

Since $E_1$ and $E_2$ are smooth and intersect transversely, the local equation

$$\alpha\left(x_3 + \frac{\beta}{2\alpha}\right)^2 - \left(\frac{\epsilon_1\epsilon_2}{4\alpha}\right) + \delta x_4^2 = 0$$

shows that the 9 singular points are nodes of $V$. □

## 2.2 Singularities

We introduce an important property of singularities of $V$.

DEFINITION 6. We say that a complete $k$-variety $X$ has $\mathbb{L}$-*rational singularities* if there is a resolution of singularities $Y \to X$ such that

$$[Y] = [X] \mod \mathbb{L}$$

in the Grothendieck ring of $k$-varieties.

Over a field of characteristic zero, $X$ has $\mathbb{L}$-rational singularities if and only if every resolution of singularities $Y \to X$ satisfies

$$[Y] = [X] \mod \mathbb{L}.$$

This follows from the result of Larsen and Lunts, and the author does not know if the assumption on the characteristic can be removed from the previous sentence.

In the remainder of this subsection, suppose that the base field is algebraically closed of characteristic zero.

PROPOSITION 7. *$V$ has $\mathbb{L}$-rational singularities.*

We construct an explicit resolution of singularities $\pi : W \to V$ with the property
$$[W] = [V] \mod \mathbb{L}.$$
Let $\widetilde{\mathbb{P}}^4$ be the blowup $\mathbb{P}^4$ along the line $L$. The blowup resolves indeterminacies of the projection

$$\begin{array}{c} \widetilde{\mathbb{P}}^4 \\ \swarrow \quad \searrow \\ \mathbb{P}^4 \dashrightarrow^{(z_0:z_1:z_2)} \mathbb{P}^2. \end{array}$$

If we write $\widetilde{V}$ for the strict transform of $V$, then the above diagram restricts to

$$\begin{array}{c} \widetilde{V} \\ \swarrow \quad \searrow \\ V \dashrightarrow^{(z_0:z_1:z_2)} \mathbb{P}^2. \end{array}$$

LEMMA 8. $\widetilde{V}$ is smooth over $L$.

*Proof.* It is enough to verify the assertion over the affine charts $\{z_3 \neq 0\} \subseteq \mathbb{P}^4$ and $\{z_4 \neq 0\} \subseteq \mathbb{P}^4$. We do this for the latter chart. Let $x_0, x_1, x_2, x_3$ be the affine coordinates
$$x_0 := \frac{z_0}{z_4}, \quad x_1 := \frac{z_1}{z_4}, \quad x_2 := \frac{z_2}{z_4}, \quad x_3 := \frac{z_3}{z_4}.$$
The blowup over this chart is covered by three affine charts in a standard way. By symmetry between $x_0$, $x_1$, $x_2$, in the defining equation of $V$, it is enough to prove the smoothness in any one of the three affine charts.

Consider the blowup chart whose induced affine coordinates are
$$y_0 := x_0, \quad y_1 := \frac{x_1}{x_0}, \quad y_2 := \frac{x_2}{x_0}, \quad y_3 := x_3.$$
The exceptional divisor of the blowup is defined by $y_0 = 0$, and the strict transform of $V$ is defined by
$$f := \alpha(1, y_1, y_2) \cdot y_3^2 + \beta(1, y_1, y_2) \cdot y_0 y_3 + \gamma(1, y_1, y_2) \cdot y_0^2 + \delta(1, y_1, y_2) = 0.$$
A singular point $p$ of $\widetilde{V}$ over $L$ in this chart satisfies
$$\frac{\partial f}{\partial y_1} = \frac{\partial \alpha}{\partial y_1} \cdot y_3^2 + \frac{\partial \delta}{\partial y_1} = 0, \quad \frac{\partial f}{\partial y_2} = \frac{\partial \alpha}{\partial y_2} \cdot y_3^2 + \frac{\partial \delta}{\partial y_2} = 0, \quad \frac{\partial f}{\partial y_3} = \alpha \cdot 2 y_3 = 0.$$

1. If $y_3(p) = 0$, then the above equations read
$$\frac{\partial \delta}{\partial y_1}(p) = 0, \quad \frac{\partial \delta}{\partial y_2}(p) = 0, \quad \text{and} \quad \delta(p) = 0.$$

This contradicts that $D$ is smooth.

2. If $y_3(p) \neq 0$, then the above equations read
$$\operatorname{rank} \begin{pmatrix} \frac{\partial \alpha}{\partial y_1}(p) & \frac{\partial \delta}{\partial y_1}(p) \\ \frac{\partial \alpha}{\partial y_2}(p) & \frac{\partial \delta}{\partial y_2}(p) \end{pmatrix} \leqslant 1 \quad \text{and} \quad \alpha(p) = \delta(p) = 0.$$

This contradicts that $A$ and $D$ intersect transversely.

Therefore $\widetilde{V}$ is smooth over $L \cap \{z_4 \neq 0\}$. The argument over $L \cap \{z_3 \neq 0\}$ is similar. □

Let $S$ be the exceptional surface of the blowup $\widetilde{V} \to V$.

**LEMMA 9.** *$S$ is smooth and rational.*

*Proof.* Write $w_0, w_1, w_2$ for the homogeneous coordinates of $\mathbb{P}^2$. By construction, $S$ is a hypersurface of $L \times \mathbb{P}^2$ defined by the bihomogeneous polynomial

$$\alpha(w_0, w_1, w_2) \cdot z_3^2 + \delta(w_0, w_1, w_2) \cdot z_4^2 = 0.$$

We show the smoothness of $S$ in the chart with affine coordinates

$$x_1 := \frac{w_1}{w_0}, \quad x_2 := \frac{w_2}{w_0}, \quad x_3 := \frac{z_3}{z_4}.$$

A singular point $p$ of $S$ in this chart satisfies

$$g := \alpha(1, x_1, x_2) \cdot x_3^2 + \delta(1, x_1, x_2) = 0$$

and its partial derivatives

$$\frac{\partial g}{\partial x_1} = \frac{\partial \alpha}{\partial x_1} \cdot x_3^2 + \frac{\partial \delta}{\partial x_1} = 0, \quad \frac{\partial g}{\partial x_2} = \frac{\partial \alpha}{\partial x_2} \cdot x_3^2 + \frac{\partial \delta}{\partial x_2} = 0, \quad \frac{\partial g}{\partial x_3} = \alpha \cdot 2x_3 = 0.$$

1. If $x_3(p) = 0$, then the above equations read

$$\frac{\partial \delta}{\partial x_1}(p) = 0, \quad \frac{\partial \delta}{\partial x_2}(p) = 0, \quad \text{and} \quad \delta(p) = 0.$$

This contradicts that $D$ is smooth.

2. If $x_3(p) \neq 0$, then the above equations read

$$\operatorname{rank} \begin{pmatrix} \frac{\partial \alpha}{\partial x_1}(p) & \frac{\partial \delta}{\partial x_1}(p) \\ \frac{\partial \alpha}{\partial x_2}(p) & \frac{\partial \delta}{\partial x_2}(p) \end{pmatrix} \leqslant 1 \quad \text{and} \quad \alpha(p) = \delta(p) = 0.$$

This contradicts that $A$ and $D$ intersect transversely.

Therefore $S$ is smooth in this affine chart. The smoothness of $S$ in other charts can be checked in a similar way.

To complete the proof, note from its defining equation that $S$ is a conic bundle over $L$. Since the base field is algebraically closed, a conic bundle over a rational curve is rational by Tsen's theorem [Lan52]. $\square$

Let $W$ be the blowup of all the nodes of $\widetilde{V}$. This defines a resolution of singularities of $V$ which fits into the commutative diagram

$$\begin{array}{ccc} & W & \\ {\scriptstyle \pi} \swarrow & & \searrow {\scriptstyle \mathrm{pr}} \\ V \dashleftarrow\text{-}\text{-}\xrightarrow{(z_0:z_1:z_2)}\text{-}\text{-}\dashrightarrow & & \mathbb{P}^2. \end{array}$$

*Proof of Proposition 7.* The first blowup $\widetilde{V} \to V$ replaces the line $L$ with the surface $S$. Since $S$ is smooth and rational (Lemma 9), we have

$$[\widetilde{V}] = [V] \mod \mathbb{L}$$

in the Grothendieck ring of varieties. The second blowup $W \to \widetilde{V}$ replaces each node by the projectivization of its tangent cone, which can be identified with a smooth quadric surface. Since any smooth quadric surface is rational over an algebraically closed base field, we have

$$[W] = [\widetilde{V}] \mod \mathbb{L}$$

in the Grothendieck ring of varieties. Therefore, $V$ has $\mathbb{L}$-rational singularities. □

### 2.3 Torsion

Let $\pi : W \to V$ be the resolution of singularities constructed in the previous subsection. We work over the field of complex numbers.

PROPOSITION 10. *We have*

(i) $\operatorname{Tors} H^4(V; \mathbb{Z}) \simeq \operatorname{Tors} H_3(V; \mathbb{Z}) \neq 0$, *and*

(ii) $\operatorname{Tors} H^4(W; \mathbb{Z}) \simeq \operatorname{Tors} H_3(W; \mathbb{Z}) \neq 0$.

*It follows from Poincaré duality that* $\operatorname{Tors} H^3(W; \mathbb{Z}) \simeq \operatorname{Tors} H_2(W; \mathbb{Z}) \neq 0$.

The projective hypersurface $V$ is one of few examples of a complete intersection in $\mathbb{P}^n$ which has torsion in its middle homology group. Other known examples are

- normal cubic surfaces with singularities of type $4A_1$, $3A_2$, $A_1A_5$, $2A_1A_3$ [Dim92], and
- certain intersections of quadric hypersurfaces [Dim86, GP08].

It is worth noting that there is an algebraic cycle whose class is a nonzero torsion element in $H_2(W; \mathbb{Z})$. This is in contrast to, for example, the toric hypersurfaces with torsion constructed in [BK06].

*Proof.* The indicated isomorphisms between homology and cohomology torsion subgroups come from the universal coefficient theorem for cohomology. We claim that (i) and (ii) are equivalent. More precisely, we claim

$$\operatorname{Tor} H_3(W; \mathbb{Z}) \simeq \operatorname{Tor} H_3(V; \mathbb{Z}).$$

Write $E \subseteq W$ for the inverse image of the singular locus $V_{\text{sing}} \subseteq V$, and consider the Borel-Moore homology exact sequences

$$\cdots \longrightarrow H_3(E; \mathbb{Z}) \longrightarrow H_3(W; \mathbb{Z}) \longrightarrow H_3^{BM}(W \setminus E; \mathbb{Z}) \longrightarrow H_2(E; \mathbb{Z}) \longrightarrow \cdots$$
$$\cdots \longrightarrow H_3(V_{\text{sing}}; \mathbb{Z}) \longrightarrow H_3(V; \mathbb{Z}) \longrightarrow H_3^{BM}(V \setminus V_{\text{sing}}; \mathbb{Z}) \longrightarrow H_2(V_{\text{sing}}; \mathbb{Z}) \longrightarrow \cdots$$

Since $E$ is a disjoint union of 10 smooth rational surfaces, we have

$$H_3(E; \mathbb{Z}) \simeq H_3(V_{\text{sing}}; \mathbb{Z}) \simeq 0 \quad \text{and} \quad \operatorname{Tor} H_2(E; \mathbb{Z}) \simeq \operatorname{Tor} H_2(V_{\text{sing}}; \mathbb{Z}) \simeq 0.$$

It follows that

$$\operatorname{Tor} H_3(W; \mathbb{Z}) \simeq \operatorname{Tor} H_3^{BM}(W \setminus E; \mathbb{Z}) \simeq \operatorname{Tor} H_3^{BM}(V \setminus V_{\text{sing}}; \mathbb{Z}) \simeq \operatorname{Tor} H_3(V; \mathbb{Z}).$$

To show (ii), we use Artin and Mumford's "brutal procedure" for constructing torsion cycles [AM72, Section 2]. The argument in [AM72] applies to our case with little change. Recall the commutative diagram

$$\begin{array}{ccc} & W & \\ \pi \swarrow & & \searrow \text{pr} \\ V & \underset{(z_0:z_1:z_2)}{\dashrightarrow} & \mathbb{P}^2 \end{array}$$

and the equation

$$h = \alpha(z_0, z_1, z_2) \cdot z_3^2 + \beta(z_0, z_1, z_2) \cdot z_3 + \gamma(z_0, z_1, z_2) + \delta(z_0, z_1, z_2) \cdot z_4^2$$

defining $V$. A computation shows that there are four types of fibers of $\mathrm{pr}: W \to \mathbb{P}^2$:

1. If $a$ is contained in none of $E_1$, $E_2$, $D$, then
$$\mathrm{pr}^{-1}(a) \simeq \mathbb{P}^1,$$
the strict transform of a line.

2. If $a$ is contained in exactly one of $E_1$, $E_2$, $D$, then
$$\mathrm{pr}^{-1}(a) \simeq \mathbb{P}^1 \vee \mathbb{P}^1,$$
the strict transform of two lines intersecting transversely at one point.

3. If $a$ is contained in $E_1 \cap D$ or $E_2 \cap D$, then
$$\mathrm{pr}^{-1}(a) \simeq 2\mathbb{P}^1,$$
the strict transform of a double line.

4. If $a$ is contained in $E_1 \cap E_2$, then
$$\mathrm{pr}^{-1}(a) \simeq \mathbb{P}^1 \vee S(a) \vee \mathbb{P}^1,$$
where $S(a)$ is the exceptional surface of the node corresponding to $a$.

Let $i = 1, 2$. Over a point $a \in E_i \backslash A$, we have the factorization
$$h = \alpha \left( z_3 + \frac{\beta}{2\alpha} + i\sqrt{\frac{\delta}{\alpha}} z_4 \right) \left( z_3 + \frac{\beta}{2\alpha} - i\sqrt{\frac{\delta}{\alpha}} z_4 \right),$$
corresponding to the decomposition $\mathrm{pr}^{-1}(a) \simeq \mathbb{P}^1 \vee \mathbb{P}^1$. We name the components of $\mathrm{pr}^{-1}(a)$ and write
$$\mathrm{pr}^{-1}(a) =: L_1(a) \vee L_2(a).$$
Note that the two components are indexed by the square roots $\pm\sqrt{\frac{\delta}{\alpha}}$. Over points in $E_i \cap A$, the above factorization specializes to
$$h = \left(\sqrt{\gamma} + i\sqrt{\delta} z_4\right)\left(\sqrt{\gamma} - i\sqrt{\delta} z_4\right).$$
Consider the double covering of $E_i$ obtained by adjoining $\sqrt{\frac{\delta}{\alpha}}$ to the function field of $E_i$. This double covering, ramified over $E_i \cap D$, parametrizes the two components of $\mathrm{pr}^{-1}(a)$ over $E_i$.

We choose a base point $a_i \in E_i \backslash D$. It is not difficult to see that there are oriented loops
$$\sigma_i' \simeq S^1 \quad \text{and} \quad \sigma_i'' \simeq S^1,$$
containing $a_i$ and contained in $E_i \backslash D$, such that

- moving $\mathrm{pr}^{-1}(a)$ around $\sigma_i'$, $L_1(a)$ and $L_2(a)$ are interchanged,
- moving $\mathrm{pr}^{-1}(a)$ around $\sigma_i''$, $L_1(a)$ and $L_2(a)$ are not interchanged, and
- $\sigma_i'$ and $\sigma_i''$ intersect only at $a_i$, and this intersection is transverse in $E_i$.

In fact, without loss of generality, we may picture $\sigma_i'$ and $\sigma_i''$ as the meridian and longitude of a real torus under some identification $E_i \simeq S^1 \times S^1$. Furthermore, we may suppose that $\sigma_i'$ and $\sigma_i''$ do not meet $E_1 \cap E_2$.

Define three oriented cycles
$$\xi_3 := \bigcup_{a \in \sigma_1''} L_2(a), \quad \widetilde{\xi}_3 := \bigcup_{a \in \sigma_1''} L_1(a), \quad \text{and} \quad \xi_2 := L_1(a_1) - L_1(a_2).$$

The first and the second are cycles because $\sigma_1''$ does not interchange $L_1(a)$ and $L_2(a)$. The orientation of $\xi_2$ is the canonical one, and the orientations of $\xi_3$ and $\widetilde{\xi}_3$ are given by the product

$$(\text{orientation of } \sigma_1'') \times (\text{canonical orientations of } L_i(a)).$$

We claim that $2\xi_2$ is homologous to zero in $W$. To see this, choose a closed oriented interval $\sigma_{12}'$ in $\mathbb{P}^2$ joining a point of $\sigma_1'$ with a point of $\sigma_2'$, whose interior does not intersect $E_1$, $E_2$, and $D$. There is a family of oriented 2-dimensional cycles in $W$ over

$$\sigma_1' \cup \sigma_{12}' \cup \sigma_2',$$

connecting $2L_1(a_1)$ with $2L_1(a_2)$ as oriented cycles. The orientations on $\sigma_1'$, $\sigma_{12}'$, $\sigma_2'$, together with canonical orientations of the members of the family, can be used to make the total space of this family an oriented chain whose boundary is $2\xi_2$. We denote this 3-dimensional oriented chain by $\eta_3$.

Next we claim that $2\xi_3$ is homologous to zero in $W$. Moving $\xi_3$ by dragging $\sigma_1''$ along $\sigma_1'$ in $E_1$, we see that $\xi_3$ and $\widetilde{\xi}_3$ are homologous to each other. It follows that $2\xi_3$ is homologous to the inverse image

$$\mathrm{pr}^{-1}(\sigma_1'') = \xi_3 + \widetilde{\xi}_3,$$

which is homologous to zero in $W$ because $\sigma_1''$ is homologous to zero in $\mathbb{P}^2$.

Finally, we show that $\xi_2$ and $\xi_3$ are not homologous to zero in $W$, using the *torsion linking form*

$$\mathscr{L} : \mathrm{Tor}\, H_2(W;\mathbb{Z}) \times \mathrm{Tor}\, H_3(W;\mathbb{Z}) \longrightarrow \mathbb{Q}/\mathbb{Z}.$$

See [ST80, Section 77]. By definition, the value of $\mathscr{L}$ at $(\xi_2, \xi_3)$ can be computed using any oriented chain whose boundary is a multiple of $\xi_2$. Using $\eta_3$ in particular to compute the value of $\mathscr{L}$, we have

$$\mathscr{L}(\xi_2, \xi_3) = \frac{\eta_3 \cdot \xi_3}{2}.$$

Since $\sigma_1'$ and $\sigma_1''$ intersect transversely in $E_1$ at one point, $\eta_3$ and $\xi_3$ meet transversely in $W$ at one point, the singular point of $\mathrm{pr}^{-1}(a_1)$. It follows that $\mathscr{L}(\xi_2, \xi_3)$ is nonzero in $\mathbb{Q}/\mathbb{Z}$, and hence $\xi_2$ and $\xi_3$ are not homologous to zero in $W$. □

COROLLARY 11. *$W$ and $V$ are not stably rational.*

*Proof.* It is shown in [AM72, Proposition 1] that the torsion subgroup of $H^3(X;\mathbb{Z})$ is a birational invariant of a smooth complete variety $X$. The Künneth theorem for cohomology shows that the torsion subgroup is in fact a stable birational invariant of smooth complete varieties. It follows from Proposition 10 that $W$ is not stably rational. □

*Remark* 12. Following [CM87], we define a *Fano threefold* to be a three dimensional embedded projective variety whose general linear curve section is canonically embedded. When the variety is smooth, this definition, originally used by Fano, agrees with the condition that the anticanonical class is very ample.

Conte and Murre show in [CM87, Proposition 2.2] that a projective threefold $X \subseteq \mathbb{P}^n$ is a Fano threefold if and only if

(i) $X$ is Cohen-Macaulay,
(ii) $X$ is projectively normal, and
(iii) the dualizing sheaf $\omega_X$ is isomorphic to $\mathcal{O}_X(-1)$.

We note that $V \subseteq \mathbb{P}^4$ is a Fano threefold. However, as Iskovskikh remarks in [Isk84, Isk97], a smooth Fano threefold $X$ has no torsion in $H^3(X; \mathbb{Z})$. It follows that the Fano threefold $V$ is not stably birational to any smooth Fano threefold.

*Remark* 13. A variety is said to be *rationally connected* if there is a rational curve in the variety joining any two general points. The quartic threefold $V$ is rationally connected because it is a conic bundle over a rational surface. See [GHS03, Corollary 1.3].

### 2.4 Proof of Theorem 3

We prove the assertion when the base field is $\mathbb{C}$. By the theorem of Larsen and Lunts, a smooth and complete variety $X$ is stably rational if and only if

$$[X] = 1 \mod \mathbb{L}$$

in the Grothendieck ring of varieties. Since $V$ has $\mathbb{L}$-rational singularities (Proposition 7), it is enough to show that a smooth model of $V$ is not stably rational. This is the content of Corollary 11.

Next we prove the general case. Let $k$ be a field of characteristic zero. We may assume that $k$ is algebraically closed because there is a natural ring homomorphism

$$K_0(\mathrm{Var}_k)/(\mathbb{L}) \longrightarrow K_0(\mathrm{Var}_{\overline{k}})/(\mathbb{L}).$$

We deduce the assertion for $\overline{k}$ from that for $\mathbb{C}$ using the Lefschetz principle. One can check that the equivalent statement on stable rationality of $V$ can be expressed in a countable conjunction of first-order sentences. This is enough to deduce the general case from the special case using the completeness of $\mathrm{ACF}_0$ [Mar02, Corollary 3.2.3].

## References

AM72    M. Artin and D. Mumford, *Some elementary examples of unirational varieties which are not rational*, Proc. London Math. Soc. (3) **25** (1972), 75–95.

Ax64    J. Ax, *Zeroes of polynomials over finite fields*, Amer. J. Math. **86** (1964), 255–261.

BK06    V. Batyrev and M. Kreuzer, *Integral cohomology and mirror symmetry for Calabi-Yau 3-folds*, Mirror symmetry V, 255–270, AMS/IP Stud. Adv. Math. **38**, Amer. Math. Soc., Providence, RI, 2006.

BCSS85    A. Beauville, J.-L. Colliot-Thélène, J.-J. Sansuc, and P. Swinnerton-Dyer, *Variétés stablement rationnelles non rationnelles*, Ann. of Math. (2) **121** (1985), 283–318.

Bil13    E. Bilgin, *On the Classes of Hypersurfaces of Low Degree in the Grothendieck Ring of Varieties*, Int. Math. Res. Not. IMRN, to appear.

Bil12    E. Bilgin, *Classes of Some Hypersurfaces in the Grothendieck Ring of Varieties*, PhD Thesis, Universität Duisburg-Essen, 2012, http://d-nb.info/1029288461.

BS12    F. Brown and O. Schnetz, *A K3 in $\phi^4$*, Duke Math. J. **161** (2012), 1817–1862.

Cay69a    A. Cayley, *A Memoir on Quartic Surfaces*, Proc. London Math. Soc. (1869), 19–69.

Cay69b    A. Cayley, *Second Memoir on Quartic Surfaces*, Proc. London Math. Soc. (1869), 198–202.

Cay69c    A. Cayley, *Third Memoir on Quartic Surfaces*, Proc. London Math. Soc. (1869), 233–266

CM87    A. Conte and J. P. Murre, *On the definition and on the nature of the singularities of Fano threefolds*, Conference on algebraic varieties of small dimension, Rend. Sem. Mat. Univ. Politec. Torino 1986, Special Issue, 51–67, 1987.

Cos83    F. Cossec, *Reye congruences*, Trans. Amer. Math. Soc. **280** (1983), 737–751.

# Geometric Chevalley-Warning conjecture

Dim86     A. Dimca, *On the homology and cohomology of complete intersections with isolated singularities*, Compositio Math. **58** (1986), 321–339.

Dim92     A. Dimca, *Singularities and Topology of Hypersurfaces*, Universitext, Springer-Verlag, New York, 1992.

GHS03     T. Graber, J. Harris, and J. Starr, *Families of rationally connected varieties*, J. Amer. Math. Soc. **16** (2003), 57–67.

GP08     M. Gross and S. Pavanelli, *A Calabi-Yau threefold with Brauer group* $(\mathbb{Z}/8\mathbb{Z})^2$, Proc. Amer. Math. Soc. **136** (2008), 1–9.

Isk84     V. A. Iskovskih, *Algebraic threefolds with special regard to the problem of rationality*, Proceedings of the International Congress of Mathematicians (Warsaw, 1983), 733–746, PWN, Warsaw, 1984.

Isk97     V. A. Iskovskikh, *On the rationality problem for three-dimensional algebraic varieties*, Proceedings of the Steklov Mathematics Institute **218** (1997), 186–227.

Jes16     C. M. Jessop, *Quartic surfaces with singular points*, Cambridge University Press, Cambridge, 1916.

Kar     I. Karzhemanov, *On one stable birational invariant*, `arXiv:1307.5605`.

Lan52     S. Lang, *On quasi algebraic closure*, Ann. Math. **55** (1952), 373–390.

LL03     M. Larsen and V. A. Lunts, *Motivic measures and stable birational geometry*, Mosc. Math. J. **3** (2003), 85–95.

Laz04     R. Lazarsfeld, *Positivity in Algebraic Geometry I*, Ergebnisse der Mathematik und ihrer Grenzgebiete. 3. Folge. A Series of Modern Surveys in Mathematics **48**, Springer-Verlag, Berlin, 2004.

Lia13     X. Liao, *Stable birational equivalence and geometric Chevalley-Warning*, Proc. Amer. Math. Soc., to appear.

Mar02     D. Marker, *Model Theory. An Introduction*, Graduate Texts in Mathematics **217**, Springer-Verlag, New York, 2002.

Ngu12     L. D. T. Nguyen, *Unramified cohomology, $\mathbb{A}^1$-connectedness, and the Chevalley-Warning problem in Grothendieck ring*, C. R. Acad. Sci. Paris, Ser. I **350** (2012), 613–615.

ST80     H. Seifert and W. Threlfall, *A Textbook of Topology*, Pure and Applied Mathematics **89**. Academic Press, Inc., New York-London, 1980.

June Huh     junehuh@umich.edu

Department of Mathematics, University of Michigan, Ann Arbor, MI 48109, USA